# Reducible Conformal Minimal Immersion with Constant Curvature from $S^2$ to $Q_6$

**Jiao Xiaoxiang[1]. Li Mingyue[2]**

School of Mathematical Science, University of Chinese Academy of Sciences, Beijing 100049, China

**Abstract** The geometry of conformal minimal two-spheres immersed in $G(2,6;\mathbb{R})$ is studied in this paper by harmonic maps. Then in most cases, we determine the linearly full reducible conformal minimal immersions from $S^2$ to $G(2,8;\mathbb{R})$ identified with the complex hyperquadric $Q_6$. We also give some examples, up to an isometry of $G(2,8;\mathbb{R})$, in which none of the spheres are congruent, with the same Gaussian curvature.

## 1 Introduction

It is a meaningful problem in studying the differential geometry of conformal minimal two-spheres with constant curvature in different Riemannian symmetric spaces, which also has a long history. In 1988, Bolton et al.[1] studied the conformal minimal immersions with constant curvature from a two-sphere $S^2$ to a complex projective space $\mathbb{C}P^n$, and proved the rigid theorem that all linearly full conformal minimal immersions with constant curvature from $S^2$ to $\mathbb{C}P^n$, up to an isometry of $\mathbb{C}P^n$, belongs to the Veronese sequence. However, the similar rigidity problems in general Riemannian symmetric space are not easy. The typical examples are the complex Grassmann manifold $G(k,n;\mathbb{C})$ and complex hyperquadric $Q_n$. In spite of this, there are many important works which have been done for these years. Li and Yu[12] found all the minimal two-spheres with constant curvature in $G(2,4;\mathbb{C})$. Wang and Jiao[15] considered a conformal minimal immersion from $S^2$ into $Q_2$, and proved that its Gaussian curvature $K$ and normal curvature $K^\perp$ satisfy $K+K^\perp=4$ and that the ellipse of curvature is a circle. And then Li et al.[11] gave a classification theorem of linearly full totally unramified conformal minimal immersions with constant curvature from $S^2$ to $Q_3$. They also pointed that all such immersions can be presented by Veronese sequence in $\mathbb{C}P^4$. Furthermore, for general linearly full totally unramified conformal minimal two-spheres immersed in complex hyperqudric $Q_n$, Jiao and Li [10] discussed all reducible cases of isotropy order $r \geq n-2$ and all irreducible cases of finite isotropy order $r \geq n-3$. Recently, Jiao and Li M.Y[13] established a classification theorem of linearly full conformal minimal two-spheres with constant curvature in $Q_4$. Jiao and Li Hong[14] determined the Gaussian curvature of all linearly full totally unramified irreducible and all linearly full reducible conformal minimal immersions with constant curvature from $S^2$ to $Q_5$.

It is well known that $Q_{n-2}$ can be identified with $G(2,n;\mathbb{R})$ in $\mathbb{C}P^{n-1}$, which is considered as a totally geodesic submanifold in $G(2,n;\mathbb{C})$. In 1986, Burstall and Wood[3] gave the explicit construction of all two-spheres in $G(2,n;\mathbb{C})$ and proved that any harmonic maps from $S^2$ to $G(2,n;\mathbb{C})$ can be obtained by a harmonic map, a Frenet pair or a mixed pair. Then in 1989 Baly-El-Dien and Wood[2] constructed all the harmonic two-spheres with finite isotropy order in $G(2,n;\mathbb{R})$. The purpose of this paper is to discuss the linearly full reducible harmonic maps from $S^2$ to $G(2,8;\mathbb{R})$ with constant curvature by using the theory of harmonic maps they gave.

The arrangement of this paper is as follows. In Section 2, at first, we identify $Q_{n-2}$ with $G(2,n;\mathbb{R})$, and then from the view of harmonic sequences, we state some fundamental formulas and results in terms of $G(k,n;\mathbb{C})$. At last, we review some properties of Veronese sequence and the rigidity theorem in $\mathbb{C}P^n$. In Section 3, we apply Bahy-El-Dien and Wood's[2] method to discuss some geometric properties of reducible harmonic maps from $S^2$ into $G(2,8;\mathbb{R})$ with constant curvature, and give a classification theorem of linearly full conformal minimal two-spheres in $G(2,8;\mathbb{R})$ under some conditions (see Theorem 1 below).

## 2 Preliminaries

(A) For $0<k<n$, we consider complex Grassmann manifold $G(k,n;\mathbb{C})$ as the set of Hermitian orthogonal projections from $\mathbb{C}^n$ onto a $k$-dimensional subspace in $\mathbb{C}^n$. Here $\mathbb{C}^n$ is endowed with the Hermitian inner product defined by

$$\langle \boldsymbol{x},\boldsymbol{y}\rangle = x_1\bar{y}_1+\cdots+x_n\bar{y}_n,$$

Where $\boldsymbol{x}=(x_1,\cdots,x_n)^T$ and $\boldsymbol{y}=(y_1,\cdots,y_n)^T$ are two elements of $\mathbb{C}^n$.

Let $G(k,n;\mathbb{R})$ denote the Grassmannian of all real $k$-dimensional subspaces of $\mathbb{R}^n$ and

$$\sigma: G(k,n;\mathbb{C})\to G(k,n;\mathbb{C}),$$

denote the complex conjugation of $G(k,n;\mathbb{C})$. It is easy to verify that $\sigma$ is an isometry with the respect to the standard Riemannian metric of $G(k,n;\mathbb{C})$. $G(k,n;\mathbb{R})$ is the fixed point set of it. Thus $G(k,n;\mathbb{R})$ is a totally geodesic submanifold of $G(k,n;\mathbb{C})$.

Map

$$Q_{n-2}\to G(2,n;\mathbb{R})$$

by

$$q\to\frac{\sqrt{-1}}{2}Z\wedge\bar{Z},$$

where $q\in Q_{n-2}$ and $Z$ is a homogeneous coordinate vector of $q$. It is obvious that the map is one-to-one and onto, and it is an isometry. Thus we can identify $Q_{n-2}$ with $G(2,n;\mathbb{R})$. Here, suppose that the metric on $G(2,n;\mathbb{R})$ is given by Section 2 of [8], then the metric is twice as much as the standard metric on $Q_{n-2}$ induced by the inclusion $\tau: Q_{n-2}\to\mathbb{C}P^{n-1}$, where this latter space is endowed with the Fubini-Study metric of constant holomorphic sectional curvature 4.

(B) Next, we will introduce general expressions of some geometric qualities about conformal minimal immersions from $S^2$ to complex Grassiann manifold $G(k,n;\mathbb{C})$.

Let $U(n)$ be the unitary group, $M$ be a simply connected domain in the unit sphere $S^2$, and let $(z,\bar{z})$ be a complex coordinate on $M$. Then we take the metric $ds_M^2 = dzd\bar{z}$ on $M$. Denote

$$A_z = \frac{1}{2}t^{-1}\partial t,\ A_{\bar{z}} = \frac{1}{2}t^{-1}\bar{\partial} t,$$

where $t: M \to U(n)$ is a smooth map, $\partial = \partial/\partial z$, $\bar{\partial} = \partial/\partial \bar{z}$.

Suppose that $t: S^2 \to U(n)$ is an isometric immersion, then $t$ is conformal and minimal if it is harmonic. Let $\omega$ be the Maurer-Cartan form on $U(n)$, and $ds_{U(n)}^2 = 1/8 tr\omega\omega^*$. Then the metric induced by $t$ on $S^2$ is locally given by

$$ds^2 = -trA_z A_{\bar{z}} dzd\bar{z}.$$

Consider the complex Grassmann manifold $G(k,n;\mathbb{C})$ as the set of Hermitian orthogonal projections from $\mathbb{C}^n$ onto a $k$-dimensional subspace in $\mathbb{C}^n$. Then a map $\psi: M \to G(k,n;\mathbb{C})$ is a Hermitian orthogonal projection onto a $k$-dimensional subbundle $\underline{\psi}$ of the trivial bundle $\underline{\mathbb{C}^n} = M \times \mathbb{C}^n$ given by setting fiber $\underline{\psi}_x = \psi(x)$ for all $x \in M$. $\underline{\psi}$ is called (a) *harmonic ((sub-)bundle)* wherever $\psi$ is a harmonic map. And here $t = \psi - \psi^{\perp}$ is a map from $S^2$ into $U(n)$. It is well known that $\psi$ is harmonic if and only if $t$ is harmonic.

Let $\psi: S^2 \to G(k,n;\mathbb{C})$ be a harmonic map. Then from $\psi$, two harmonic sequences are derived as follows:

$$\underline{\psi} = \underline{\psi}_0 \xrightarrow{\partial'} \underline{\psi}_1 \xrightarrow{\partial'} \cdots \xrightarrow{\partial'} \underline{\psi}_i \xrightarrow{\partial'} \cdots, \tag{0.1}$$

$$\underline{\psi} = \underline{\psi}_0 \xrightarrow{\partial''} \underline{\psi}_{-1} \xrightarrow{\partial''} \cdots \xrightarrow{\partial''} \underline{\psi}_{-i} \xrightarrow{\partial''} \cdots, \tag{0.2}$$

where $\underline{\psi}_i = \partial' \underline{\psi}_{i-1}$ and $\underline{\psi}_{-i} = \partial'' \underline{\psi}_{-i+1}$ are Hermitian orthogonal projections from $S^2 \times \mathbb{C}^n$ onto $\underline{Im}(\psi_{i-1}^{\perp}\partial\psi_{i-1})$ and $\underline{Im}(\psi_{-i+1}^{\perp}\bar{\partial}\psi_{-i+1})$, respectively. In the following, we also denote them by $\partial^{(i)}\underline{\psi}$ and $\partial^{(-i)}\underline{\psi}$, respectively, $i = 1, 2, \cdots$.

Now recall[4] that a harmonic map $\psi: S^2 \to G(k,n;\mathbb{C})$ in (2.1) [resp.(2.2)] is said to be $\partial'$-*irreducible* (resp. $\partial''$-*irreducible*) if rank $\underline{\psi}$ = rank $\underline{\psi}_1$ (resp. rank $\underline{\psi}$ = rank $\underline{\psi}_{-1}$) and $\partial'$- *reducible* (resp. $\partial''$-*reducible*) otherwise. In particular, let $\psi$ be a harmonic map from $S^2$ to $G(2,n;\mathbb{R})$, then $\psi$ is $\partial'$-irreducible (resp. $\partial'$-reducible ) if and only if $\psi$ is $\partial''$-irreducible (resp. $\partial''$-reducible). In this case, we simply say that $\psi$ is irreducible (resp. reducible). According to Wood[6], call a harmonic map $\psi: S^2 \to G(k,n;\mathbb{C})$ (*strongly*) *isotropic* if $\psi_i \perp \psi$, $\forall i \in \mathbb{Z}$, $i \neq 0$.

**Definition 2.1** Let $\psi: S^2 \to G(k,n;\mathbb{C})$ be a map.

- If $\underline{\psi}$ cannot be contained in any proper trivial subbundle $S^2 \times \mathbb{C}^m$ of $S^2 \times \mathbb{C}^n$ $(m < n)$, $\psi$ is linearly full.
- If $\psi$ is harmonic, define its isotropy order to be the greatest integer $r$ such that $\psi_i \perp \psi$ for all $i$ with $1 \leq i \leq r$, and if $\psi_i \perp \psi$, $(\forall i \in \mathbb{Z}, i \neq 0)$, $\psi$ is (strongly) isotropic and set $r = \infty$.

In this paper, we always assume that $\psi$ is linearly full.

**Lemma 2.2**[2] *Let* $\psi: S^2 \to G(k,n;\mathbb{R})$ *be a harmonic map of finite isotropy order* $r$. *Then* $r$ *is odd. Further, writing* $r = 2s+1$, *we have* $0 \leq s \leq [(n-3)/2]$. (*Here* $[x]$ *denotes integer part of* $x$.)

Suppose that $\psi: S^2 \to G(k,n;\mathbb{C})$ is a linearly full harmonic map and it belongs to the following harmonic sequence:

$$\underline{\psi}_0 \xrightarrow{\partial'} \cdots \underline{\psi} = \underline{\psi}_i \xrightarrow{\partial'} \underline{\psi}_{i+1} \xrightarrow{\partial'} \cdots \xrightarrow{\partial'} \underline{\psi}_{i_0} \xrightarrow{\partial'} 0 \tag{0.3}$$

for some $i = 0, \cdots, i_0$. We choose the local unit orthogonal frame $e_1^{(i)}, e_2^{(i)}, \cdots, e_{k_i}^{(i)}$ such that they locally span subbundle $\underline{\psi}_i$ of $S^2 \times \mathbb{C}^n$, where $k_i = \text{rank}\, \underline{\psi}_i$.

Let $W_i = (e_1^{(i)}, e_2^{(i)}, \cdots, e_{k_i}^{(i)})$ be a $(n \times k_i)$-matrix. Then we have

$$\begin{aligned} \psi_i &= W_i W_i^*, \\ W_i^* W_i = I_{k_i \times k_i}, \quad & W_i^* W_{i+1} = 0, \quad W_i^* W_{i-1} = 0. \end{aligned} \tag{0.4}$$

By (0.4), a straightforward computation shows that

$$\begin{cases} \partial W_i = W_{i+1} \Omega_i + W_i \Phi_i, \\ \bar{\partial} W_i = -W_{i-1} \Omega_{i-1}^* - W_i \Phi_i^*, \end{cases} \tag{0.5}$$

where $\Omega_i$ is a $(k_{i+1}\times k_i)$ -matrix, $\Phi_i$ is a $(k_i\times k_i)$ -matrix for $i=0,1,2,\cdots,\ i_0$ and $\Omega_{i_0}=0$. It is evident that integrability conditions for (0.5) are

$$\bar{\partial}\Omega_i=\Phi^*_{i+1}\Omega_i-\Omega_i\Phi^*_i,$$
$$\bar{\partial}\Phi_i+\partial\Phi^*_i=\Omega^*_i\Omega_i+\Phi^*_i\Phi_i-\Omega_{i-1}\Omega^*_{i-1}-\Phi_i\Phi^*_i.$$

Set $L_i=tr(\Omega_i\Omega^*_i)$, the metric induced by $\psi_i$ is given in the form

$$ds^2_i=(L_{i-1}+L_i)\,dzd\bar{z}\,. \tag{0.6}$$

From (0.6), suppose $ds^2=\lambda^2dzd\bar{z}$ is the induced metric by $\psi$, $K$ and $B$ are its Gauss curvature and second fundamental form, respectively. Then we have

$$\begin{cases}\lambda^2=tr\,\partial\psi\,\bar{\partial}\psi,\\ K=-\dfrac{2}{\lambda^2}\partial\bar{\partial}\log\lambda^2,\\ \parallel B\parallel^2=4trPP^*\end{cases} \tag{0.7}$$

where $P=\partial(A_z/\lambda^2)$ with $A_z=(2\psi-I)\partial\psi$, $A_{\bar{z}}=(2\psi-I)\bar{\partial}\psi$, and $I$ is the identity matrix, then $P^*=-\bar{\partial}(A_{\bar{z}}/\lambda^2)$.

(C) In this section, we state the definition of degree of a smooth map from a compact Riemann surface into a complex projective space and then review the rigidity theorem of conformal minimal immersion of constant curvature from $S^2$ to $\mathbb{C}P^n$.

Suppose that $\phi:S^2\to\mathbb{C}P^n$ is a linearly full conformal minimal immersion. Then the following harmonic sequence in $\mathbb{C}P^n$ is uniquely determined by $\phi$:

$$0\xrightarrow{\partial'}\underline{\phi}_0^{(n)}\xrightarrow{\partial'}\cdots\underline{\phi}=\underline{\phi}_i^{(n)}\xrightarrow{\partial'}\cdots\xrightarrow{\partial'}\underline{\phi}_n^{(n)}\xrightarrow{\partial'}0 \tag{0.8}$$

for some $i=0,1,\cdots,n$.

Define a sequence $f_0^{(n)},\cdots,f_n^{(n)}$ of local sections of $\underline{\phi}_0^{(n)},\cdots,\underline{\phi}_n^{(n)}$ inductively so that $f_0^{(n)}$ is a nowhere zero local section of $\underline{\phi}_0^{(n)}$ (we assume without loss of generality that $\bar{\partial}f_0^{(n)}\equiv 0$) and

$$f_{(i+1)}^{(n)}=\pi_{\underline{\phi}_i^{(n)\perp}}(\partial f_i^{(n)}),\qquad i=0,1,\cdots,n-1.$$

Then we have some formulas as follows:

$$\partial f_i^{(n)}=f_{i+1}^{(n)}+\frac{\langle\partial f_i^{(n)},f_i^{(n)}\rangle}{|f_i^{(n)}|^2}f_i^{(n)},\quad i=0,\cdots,n,$$

$$\bar{\partial} f_i^{(n)} = -\frac{|f_i^{(n)}|^2}{|f_{i-1}^{(n)}|^2} f_{i-1}^{(n)}, \quad i = 1, \cdots, n.$$

Define $l_i^{(n)}$ by

$$l_i^{(n)} = \frac{|f_{i+1}^{(n)}|^2}{|f_i^{(n)}|^2}, \quad i = 0, 1, \cdots, n-1, \qquad l_{-1}^{(n)} = l_n^{(n)} = 0.$$

Then Bolton et al[1] derived the unintegrated Plücker formulae for the $l_i$'s, namely

$$\partial \bar{\partial} \log l_i = l_{i+1} - 2l_i + l_{i-1}, \quad i = 0, 1, \cdots, n-1.$$

Recall the harmonic sequence in (0.8), let $r(\partial') =$ sum of the indices of the singularities of $\partial'$, which is called the *ramification index* of $\partial'$ by Bolton et al[1]. Note that if $r(\partial') = 0$ in(0.8) for all $\partial'$, the harmonic sequence (0.8) is defined to be *totally unramified* [1].

**Definition 2.3** The degree $\sigma_i^{(n)}$ of $f_i^{(n)}$ is defined by

$$\sigma_i^{(n)} = \frac{1}{2\pi\sqrt{-1}} \int_{S^2} \partial \bar{\partial} \log |f_i^{(n)}|^2 d\bar{z} dz, \quad i = 0, 1, \cdots, n.$$

Consider the $i$-th osculating curve $G_i^{(n)} : S^2 \to \mathbb{C}P^{\binom{n+1}{i+1}-1}$ which has a local lift $\hat{G}_i^{(n)}$ into $\mathbb{C}P^{\binom{n+1}{i+1}}$ given by

$$\hat{G}_i^{(n)} = f_0^{(n)} \wedge f_1^{(n)} \wedge \cdots f_i^{(n)}, \quad i = 0, 1, \cdots, n.$$

We write $\hat{G}_i^{(n)} = g(z) \tilde{G}_i^{(n)}$, where $g(z)$ is the greatest common divisor of the $\binom{n+1}{i+1}$ components of $\hat{G}_i^{(n)}$. Then $\tilde{G}_i^{(n)}$ is a nowhere zero holomorphic curve, and the degree $\delta_i^{(n)}$ of $\hat{G}_i^{(n)}$ is given by

$$\delta_i^{(n)} = \frac{1}{2\pi\sqrt{-1}} \int_{S^2} \partial \bar{\partial} \log |\hat{G}_i^{(n)}|^2 d\bar{z} dz,$$

which is equal to the degree of the polynomial function $\tilde{G}_i^{(n)}$. After a simple calculation, we have

$$\delta_i^{(n)} = \frac{1}{2\pi\sqrt{-1}} \int_{S^2} l_1^{(n)} d\bar{z} dz.$$

Moreover, if (2.10) is a totally unramified harmonic sequence, then Bolton et al proved[1]

$$\delta_i^{(n)} = (i+1)(n-i). \tag{0.9}$$

In passing, we state the rigidity theorem of conformal minimal immersions of $S^2$ into $\mathbb{C}P^n$ with constant curvature as follows. Consider the Veronese sequence

$$0 \to \underline{V}_0^{(n)} \xrightarrow{\partial'} \underline{V}_1^{(n)} \xrightarrow{\partial'} \cdots \xrightarrow{\partial'} \underline{V}_n^{(n)} \xrightarrow{\partial'} 0.$$

For each $i=0,\cdots,n$, $V_i^{(n)}:S^2\to\mathbb{C}P^n$ is given by $V_i^{(n)}=(v_{i,0},\cdots,v_{i,n})^{\mathrm{T}}$, where , for $z\in S^2$ and $j=0,\cdots,n$,

$$v_{i,j}(z)=\frac{i!}{(1+z\bar{z})^i}\sqrt{\binom{n}{j}}z^{j-i}\sum_k(-1)^k\binom{j}{i-k}\binom{n-j}{k}(z\bar{z})^k, \tag{0.10}$$

$$|V_i^{(n)}|^2=\frac{n!\,i!}{(n-i)!}(1+z\bar{z})^{n-2q}.$$

Each map $V_i^{(n)}$ is a conformal minimal immersion with the induced metric

$$ds_i^2=\frac{n+2i(n-i)}{(1+z\bar{z})^2}dzd\bar{z},$$

and the corresponding constant curvature $K_i$ is given by

$$K_i=\frac{4}{n+2i(n-i)}.$$

By Calabi's rigidity theorem, Bolton et al. proved that the following rigidity result[1].

**Lemma 2.4** *Let $\phi:S^2\to\mathbb{C}P^n$ be a linearly full conformal minimal immersion of constant curvature. Then, up to a holomorphic isometry of $\mathbb{C}P^n$, the harmonic sequence determined by $\phi$ is Veronese sequence.*

# 3 Reducible harmonic maps of constant curvature from $S^2$ to $G(2,n;\mathbb{R})$

In the section, we see harmonic maps from $S^2$ to $G(2,8;\mathbb{R})$ as conformal minimal immersions of $S^2$ in $G(2,8;\mathbb{R})$. And then we discuss the harmonic maps only consider by reducible case.

Let $\psi:S^2\to G(2,8;\mathbb{R})$ is a linearly full reducible harmonic map with constant curvature. From Bahy-El-Dien and Wood[2], there are two cases that we should distinguish for consideration:

(1) $\psi$ is a real mixed pair with finite isotropy order;

(2) $\psi$ is (strongly) isotropic.

Firstly, we discuss the first case when $\psi$ has finite isotropy order. And suppose the finite isotropy order of $\psi$ is $r$. It follows from [12] that all linearly full reducible harmonic maps from $S^2$ to $G(2,8;\mathbb{R})$ with constant curvature of isotropy order $r\geq 4$ have been characterized. And from Lemma 2.1, we only need to consider $r=1$ and $r=3$. Then $\psi$ can be characterized by harmonic maps to $\mathbb{C}P^m (m \le 7)$. In fact,

$$\underline{\psi} = \overline{\underline{f}}_0^{(m)} \oplus \underline{f}_0^{(m)},$$

where $f_0^{(m)}: S^2 \to \mathbb{C}P^m$ is holomorphic. By using $\psi$, a harmonic sequence is derived as follows:

$$0 \xleftarrow{\partial''} \overline{\underline{f}}_m^{(m)} \xleftarrow{\partial''} \cdots \xleftarrow{\partial''} \overline{\underline{f}}_1^{(m)} \xleftarrow{\partial''} \underline{\psi} \xrightarrow{\partial'} \underline{f}_1^{(m)} \xrightarrow{\partial'} \cdots \xrightarrow{\partial'} \underline{f}_m^{(m)} \xrightarrow{\partial'} 0, \tag{0.11}$$

where $0 \xrightarrow{\partial'} \underline{f}_0^{(m)} \xrightarrow{\partial'} \underline{f}_1^{(m)} \cdots \xrightarrow{\partial'} \underline{f}_m^{(m)} \xrightarrow{\partial'} 0$ is a linearly full harmonic sequence in $\mathbb{C}P^m \subset \mathbb{C}P^7$ satisfying

$$\begin{cases} \langle f_0^{(m)}, \overline{f}_i^{(m)} \rangle = 0 \quad (0 \le i \le r), \\ \langle f_0^{(m)}, \overline{f}_{r+1}^{(m)} \rangle \ne 0 \end{cases}. \tag{0.12}$$

The induced metric of $\psi$ is given by

$$ds^2 = 2 l_0^{(m)} dz d\bar{z}, \tag{0.13}$$

where $l_0^{(m)} dz d\bar{z}$ is the induced metric of $f_0^{(m)}: S^2 \to \mathbb{C}P^m$. Since $\underline{\psi}$ is of constant curvature, by (0.13) we get that the curvature $K$ of $\psi$ satisfies

$$K = \frac{2}{m}.$$

By Lemma 2.4, up to a holomorphic isometry of $\mathbb{C}P^7$, $f_0^{(m)}$ is Veronese surface. Then we can choose a proper complex coordinate $z$ on $\mathbb{C} = S^2 \setminus \{pt\}$ so that $f_0^{(m)} = U V_0^{(m)}$, where $U \in U(8)$ and $V_0^{(m)}$ is the standard Veronese expression given in Sect.2 (adding zeros to $V_0^{(m)}$ such that $V_0^{(m)} \in \mathbb{C}^8$). Then (0.12) becomes

$$\begin{cases} \langle U V_0^{(m)}, \overline{U V}_i^{(m)} \rangle = 0 \quad (0 \le i \le r), \\ \langle U V_0^{(m)}, \overline{U V}_{r+1}^{(m)} \rangle \ne 0, \end{cases}$$

which is equivalent to

$$\begin{cases} tr W V_0^{(m)} V_i^{(m)\mathrm{T}} = 0 \quad (0 \le i \le r), \\ tr W V_0^{(m)} V_{r+1}^{(m)\mathrm{T}} \ne 0, \end{cases} \tag{0.14}$$

where $W = U^{\mathrm{T}} U$. It satisfies $W \in U(8)$ and $W = W^{\mathrm{T}}$. Might as well set

$$G_W := \{ U \in U(8) | U^{\mathrm{T}} U = W \}.$$

For a given $W$, the following statement can be easily checked:

(i). For any $U \in G_W$, $A \in SO(8)$, we have $AU \in G_W$;

(ii). For any $U, V \in G_W$, there exists a $B \in SO(8)$ such that $U = BV$.

As to the second fundamental form $B$ of $\psi$, by (0.7) and a series of calculations, we have

$$\begin{cases} \partial\psi = \dfrac{1}{|f_0^{(m)}|^2}[\bar{f}_0^{(m)}(\bar{f}_1^{(m)})^* + f_1^{(m)} f_0^{(m)*}], \\ A_z = \dfrac{1}{|f_0^{(m)}|^2}[\bar{f}_0^{(m)}(\bar{f}_1^{(m)})^* - f_1^{(m)} f_0^{(m)*}], \\ P = \dfrac{1}{2|f_1^{(m)}|^2}[\bar{f}_0^{(m)}(\bar{f}_2^{(m)})^* - f_2^{(m)} f_0^{(m)*}]. \end{cases}$$

And then we can get the following straight formula

$$\| B \|^2 = 2\frac{\delta_1^{(m)}}{\delta_0^{(m)}} - 2\frac{|\langle f_0^{(m)}, \bar{f}_2^{(m)}\rangle|^2}{|f_1^{(m)}|^4}. \tag{0.15}$$

## 3.1 Finite isotropy order $r = 1$

In this subsection, we consider $\psi$ is a linearly full reducible map with constant curvature and finite isotropy order $r = 1$.

For any integers $n \geq 3$, $s \geq 0$, $H_n^s$ denote the set of all holomorphic maps $f : S^2 \to \mathbb{C}P^{n-1}$ satisfying

$$\begin{cases} \langle \partial^{(i)} f, \bar{f} \rangle = 0 \quad (0 \leq i \leq 2s+1), \\ \langle \partial^{(2s+2)} f, \bar{f} \rangle \neq 0. \end{cases}$$

Together with (0.12), we have $f_0^{(m)} \in H_{m+1}^0$. To characterize $\psi$, we review a special case of one of Bahy-Dien and Wood's results[2]:

**Lemma 3.1** *All holomorphic maps* $f_0^{(m)} : S^2 \to \mathbb{C}P^m$ *satisfying* $f_0^{(m)} \in H_{m+1}^0$ *may be constructed by the following three steps:*

(1) *Choose* $F_0(z) : \mathbb{C} \to (\mathbb{C} \cup \{\infty\})^{m-1}$ *polynomial with* $\langle F_0(z), \bar{F}_0(z) \rangle \neq 0$.

(2) *Let* $H(z)$ *be the unique rational function* $\mathbb{C} \to (\mathbb{C} \cup \{\infty\})^{m-1}$ *with* $dH(z)/dz = F_0(z)$ *for any* $z \in \mathbb{C}$ *and* $H(0) = 0$.

(3) *Define* $F_1(z):\mathbb{C}\to\mathbb{C}^{m+1}=\mathbb{C}^{m-1}\oplus\mathbb{C}\oplus\mathbb{C}$ *by* $F_1(z)=(2H(z),1-\langle H(z),\bar{H}(z)\rangle,\sqrt{-1}(1+\langle H(z),\bar{H}(z)\rangle))$. *Then* $F_1(z)$ *is a rational function and represents the holomorphic map* $f_0^{(m)}:S^2\to\mathbb{C}P^m$ *in homogeneous coordinates.*

Using this lemma, we prove the following result.

**Lemma 3.2** *Let* $\psi:S^2\to G(2,8;\mathbb{R})$ *be a linearly full reducible harmonic map with constant curvature* $K$ *of isotropy order 1. In fact,* $\underline{\psi}=\overline{\underline{f}_0^{(m)}}\oplus\underline{f_0^{(m)}}(3\le m\le 7)$ *, where* $f_0^{(m)}:S^2\to\mathbb{C}P^m$ *is holomorphic. Suppose that* $m\le 6$*, then up to an isometry of* $G(2,8;\mathbb{R})$*,* $\psi$ *belongs to one of the following four cases:*

*(1).* $\underline{\psi}=\overline{\underline{f}}_0^{(6)}\oplus\underline{f}_0^{(6)}=\overline{\underline{UV}}_0^{(6)}\oplus\underline{UV}_0^{(6)}$ *with* $K=\dfrac{1}{3}$ *for some* $U\in G_W$*, where* $W$ *has the form* (0.29) *and*

$$w_{11},\quad w_{12},\quad -\sqrt{6}\,w_{13}-\sqrt{5}\,w_{22},\quad -4\sqrt{5}\,w_{14}-5\sqrt{6}\,w_{23},\quad -w_{15}-4w_{24}-3w_{33},$$
$$-4\sqrt{5}\,w_{25}-5\sqrt{6}\,w_{34},\quad -\sqrt{6}\,w_{35}-\sqrt{5}\,w_{44},\quad w_{45},\quad w_{55},$$

*(0.16)*

*are not all zero;*

*(2).* $\underline{\psi}=\overline{\underline{f}}_0^{(5)}\oplus\underline{f}_0^{(5)}=\overline{\underline{UV}}_0^{(5)}\oplus\underline{UV}_0^{(5)}$ *with* $K=\dfrac{2}{5}$ *for some* $U\in G_W$*, where* $W$ *has the form* (0.32) *and*

$$w_{11},\quad w_{12},\quad 3\sqrt{2}\,w_{13}+4w_{22},\quad w_{14}+3w_{23},\quad 3\sqrt{2}\,w_{24}+4w_{33},\quad w_{34},\quad w_{44}, \qquad (0.17)$$

*are not all zeros;*

*(3).* $\underline{\psi}=\overline{\underline{f}}_0^{(4)}\oplus\underline{f}_0^{(4)}=\overline{\underline{UV}}_0^{(4)}\oplus\underline{UV}_0^{(4)}$ *with* $K=\dfrac{1}{2}$ *for some* $U\in G_W$*, where* $W$ *has the form* (0.34) *and*

$$w_{11},\quad w_{12},\quad w_{13}+w_{22},\quad w_{23},\quad w_{33}, \qquad (0.18)$$

*are not all zero;*

*(4).* $\underline{\psi}=\overline{\underline{f}}_0^{(3)}\oplus\underline{f}_0^{(3)}=\overline{\underline{UV}}_0^{(3)}\oplus\underline{UV}_0^{(3)}$ *with* $K=\dfrac{2}{3}$ *for some* $U\in G_W$*, where* $W$ *has the form* (0.36) *and*

$$w_{11},\quad w_{12},\quad w_{22}, \qquad (0.19)$$

*are not all zero.*

*In each of these four cases, there are many different types of $W$, and thus there exist different $U \in U(8)$ such that $UV_0^{(m)} (m=3,4,5,6)$ are linearly full in $G(2,8;\mathbb{R})$, and they are not $SO(8)$-equivalent.*

Proof : In the beginning, it is easy to verify that $\langle f_0^{(m)}, \bar{f}_0^{(m)} \rangle = 0$ and $\langle f_1^{(m)}, \bar{f}_0^{(m)} \rangle = 0$ are equivalent. And according to the discussion, $r=1$ and $2 \le m \le 6$, in the following we deal with $W$ in case $m=6,5,4,3,2$, respectively.

(1). $m=6$, $K=\frac{1}{2}$.

Let us assume that there exists a linearly full reducible harmonic map $\underline{\psi} = \underline{\bar{f}}_0^{(6)} \oplus \underline{f_0^{(6)}} : S^2 \to G(2,8;\mathbb{R})$ with constant curvature and finite isotropy order $r=1$. Then $f_0^{(6)} \in H_7^0$ and we construct it as follows.

First of all, choose $F_0(z) = C(1,z,z^2,z^3,z^4,0)^T$, where $C=(c_{ij})_{6\times 6}$ is a constant matrix with $\langle F_0(z), \bar{F}_0(z) \rangle \neq 0$. Then by the second step, we can write $H(z)$ in the form $H(z) = \left( c_{i0} z + \frac{1}{2} c_{i1} z^2 + \frac{1}{3} c_{i2} z^3 + \frac{1}{4} c_{i3} z^4 + \frac{1}{5} c_{i4} z^5 \right)_{6\times 1}$, which gives

$$\langle H(z), \bar{H}(z) \rangle = \sum_{i=0}^{5} \left( c_{i0} z + \frac{1}{2} c_{i1} z^2 + \frac{1}{3} c_{i2} z^3 + \frac{1}{4} c_{i3} z^4 + \frac{1}{5} c_{i4} z^5 \right)^2. \tag{0.20}$$

This relationship together with the third step of Lemma (3.2), we can obtain $F_1(z)$. In fact, $F_1(z)$ is the representation of a holomorphic map of $S^2$ in $\mathbb{C}P^6$, which shows that the coefficients of $z^7$, $z^8$, $z^9$ and $z^{10}$ in (0.20) all vanish. Thus it can be expressed by

$$\begin{gathered} \frac{1}{5} \sum_{i=0}^{5} c_{i1} c_{i4} + \frac{1}{6} \sum_{i=0}^{5} 5 c_{i2} c_{13} = 0, \quad \frac{2}{15} \sum_{i=0}^{5} c_{i2} c_{i4} + \frac{1}{16} \sum_{i=0}^{5} c_{13}^2 = 0, \\ \frac{1}{10} \sum_{i=0}^{5} c_{13} c_{14} = 0, \quad \frac{1}{25} \sum_{i=0}^{5} c_{i4}^2 = 0. \end{gathered} \tag{0.21}$$

For convenience, denote

$$\langle H(z), \bar{H}(z) \rangle = \sqrt{15}\, b_2 z^2 + \sqrt{20}\, b_3 z^3 + \sqrt{15}\, b_4 z^4 + \sqrt{6}\, b_5 z^5 + b_6 z^6,$$

where $b_2$, $b_3$, $b_4$, $b_5$, $b_6$ are constant. So it is easy to see that

$$F_1(z) = \begin{pmatrix} \mathbf{D}_1 \\ \mathbf{D}_2 \end{pmatrix} = f_0^{(6)} = UV_0^{(6)},$$

where the block matrix $D_1, D_2$ can be shown as follows, and all counts of $i$ start at 0.

$$D_1 = \left(2c_{i0}z + c_{i1}z^2 + \frac{2}{3}c_{i2}z^3 + \frac{1}{2}c_{i3}z^4 + \frac{2}{5}c_{i4}z^5\right)_{6\times 1},$$

$$D_2 = \left(\sqrt{(-1)^i}(1 + (-1)^{i+1} + \sqrt{15}\,b_2 z^2 + \sqrt{20}\,b_3 z^3 + \sqrt{15}\,b_4 z^4 + \sqrt{6}\,b_5 z^5 + b_6 z^6)\right)_{2\times 1}$$

Then using the method of indeterminate coefficients, and uniting U directly, we can get

$$U = \begin{pmatrix} \mathbf{0}_{6\times 1} & & & \mathbf{A}_{6\times 5} & & & \mathbf{0}_{6\times 1} & \tilde{\mathbf{U}}_{6\times 1} \\ \frac{1}{\sqrt{2}} & 0 & -\frac{b_2}{\sqrt{2}} & -\frac{b_3}{\sqrt{2}} & -\frac{b_4}{\sqrt{2}} & -\frac{b_5}{\sqrt{2}} & -\frac{b_6}{\sqrt{2}} & u_{67} \\ \frac{\sqrt{-1}}{\sqrt{2}} & 0 & \frac{\sqrt{-1}\,b_2}{\sqrt{2}} & \frac{\sqrt{-1}\,b_3}{\sqrt{2}} & \frac{\sqrt{-1}\,b_4}{\sqrt{2}} & \frac{\sqrt{-1}\,b_5}{\sqrt{2}} & \frac{\sqrt{-1}\,b_6}{\sqrt{2}} & u_{77} \end{pmatrix}, \tag{0.22}$$

$$A = \left(\frac{2c_{i0}}{\sqrt{12}} \;\; \frac{c_{i1}}{\sqrt{30}} \;\; \frac{c_{i2}}{3\sqrt{10}} \;\; \frac{c_{i3}}{2\sqrt{30}} \;\; \frac{c_{i4}}{5\sqrt{3}}\right)_{6\times 5}, \tag{0.23}$$

where $\tilde{U} = (u_{i7})_{6\times 1}$. Set $W = U^{\mathrm{T}}U \triangleq (w_{ij})_{(8\times 8)}$. It is clear that $w_{ij} = w_{ji}$ for any $0 \le i,j \le 7$. By the expression of $V_0^{(6)}$ given in Sect.2, we get $V_0^{(6)}V_0^{(6)T}$ is a polynomial in $z$ and $\bar{z}$. Using the method of indeterminate coefficients again, the relations in (0.21) are respectively equivalent to

$$3\sqrt{10}\,w_{25} + 10\sqrt{3}\,w_{34} = 0, \quad 2\sqrt{30}\,w_{35} + \frac{15}{2}w_{44} = 0, \quad w_{45} = 0, \quad w_{55} = 0. \tag{0.24}$$

Furthermore, by the properties of isotropy order $r = 1$, we have

$$\langle F_1(z), \bar{F}_1(z)\rangle = trWV_0^{(6)}V_0^{(6)\mathrm{T}} = 0.$$

Then after a series of calculations, there exist some other relations:

$$w_{00} = w_{01} = w_{56} = w_{66} = 0, \quad w_{02} = -\frac{3}{\sqrt{15}}w_{11}, \quad w_{03} = -\frac{3}{\sqrt{2}}w_{02}, \tag{0.25}$$

$$w_{04} = -2\sqrt{2}\,w_{13} - \frac{\sqrt{15}}{2}w_{22}, \quad w_{05} = -\sqrt{5}\,w_{14} - 5\sqrt{2}\,w_{23}, \tag{0.26}$$

$$w_{06} = -6w_{15} - 15w_{24} - 10w_{33}, \quad w_{36} = -\frac{3}{\sqrt{2}}w_{45}, \quad w_{46} = -\frac{3}{\sqrt{15}}w_{55}, \tag{0.27}$$

$$w_{26} = -2\sqrt{2}\,w_{35} - \frac{\sqrt{15}}{2}w_{44}, \quad w_{16} = -\sqrt{5}\,w_{25} - 5\sqrt{2}\,w_{34}. \tag{0.28}$$

Combining the formulas (0.24)-(0.28), it is a straightforward computation to derive

$$
W=\begin{pmatrix}
0 & 0 & 0 & 0 & 0 & 0 & w_{06} & w_{07}\\
 & 0 & 0 & -\frac{\sqrt{30}}{8}w_{22} & -\frac{\sqrt{30}}{3}w_{23} & w_{15} & 0 & w_{17}\\
 & & w_{22} & w_{23} & w_{24} & -\frac{\sqrt{30}}{3}w_{34} & 0 & w_{27}\\
 & & & w_{33} & w_{34} & -\frac{\sqrt{30}}{8}w_{44} & 0 & w_{37}\\
 & & & & w_{44} & 0 & 0 & w_{47}\\
 & & * & & & 0 & 0 & w_{57}\\
 & & & & & & 0 & w_{67}\\
 & & & & & & & w_{77}
\end{pmatrix}, \tag{0.29}
$$

where

$$
w_{06}=-6w_{15}-15w_{24}-10w_{33},
$$

and terms in (0.16) are not all zero.

There are many different such type of $W$, giving different $U$. In other words, we can find many different $U$ to write $\underline{\psi}=\overline{\underline{UV}}{}_{0}^{(6)}\oplus\underline{UV}{}_{0}^{(6)}$, and they are not congruent. Now we just give one example of them:

$$
W=\begin{pmatrix}
0 & 0 & 0 & 0 & 0 & 0 & 1 & 0\\
0 & 0 & 0 & 0 & 0 & 1 & 0 & 0\\
0 & 0 & 0 & 0 & -1 & 0 & 0 & 0\\
0 & 0 & 0 & \frac{4}{5} & 0 & 0 & 0 & \frac{3}{5}\\
0 & 0 & -1 & 0 & 0 & 0 & 0 & 0\\
0 & 1 & 0 & 0 & 0 & 0 & 0 & 0\\
1 & 0 & 0 & 0 & 0 & 0 & 0 & 0\\
0 & 0 & 0 & \frac{3}{5} & 0 & 0 & 0 & -\frac{4}{5}
\end{pmatrix},
U=\begin{pmatrix}
\frac{1}{\sqrt{2}} & 0 & 0 & 0 & 0 & 0 & \frac{1}{\sqrt{2}} & 0\\
\frac{\sqrt{-1}}{\sqrt{2}} & 0 & 0 & 0 & 0 & 0 & -\frac{\sqrt{-1}}{\sqrt{2}} & 0\\
0 & \frac{1}{\sqrt{2}} & 0 & 0 & 0 & \frac{1}{\sqrt{2}} & 0 & 0\\
0 & \frac{\sqrt{-1}}{\sqrt{2}} & 0 & 0 & 0 & -\frac{\sqrt{-1}}{\sqrt{2}} & 0 & 0\\
0 & 0 & -\frac{1}{\sqrt{2}} & 0 & \frac{1}{\sqrt{2}} & 0 & 0 & 0\\
0 & 0 & \frac{\sqrt{-1}}{\sqrt{2}} & 0 & \frac{\sqrt{-1}}{\sqrt{2}} & 0 & 0 & 0\\
0 & 0 & 0 & u_{63} & 0 & 0 & 0 & u_{67}\\
0 & 0 & 0 & u_{73} & 0 & 0 & 0 & u_{77}
\end{pmatrix},
$$

where $\begin{cases} u_{63}=\frac{\sqrt{15}}{5}+\frac{\sqrt{-1}}{\sqrt{30}}\\ u_{73}=\frac{\sqrt{30}}{10}-\frac{\sqrt{-1}}{\sqrt{15}}\end{cases}$, and $\begin{cases} u_{67}=\frac{1}{\sqrt{15}}-\frac{\sqrt{-30}}{10}\\ u_{77}=\frac{1}{\sqrt{30}}+\frac{\sqrt{-15}}{5}\end{cases}$.

In this case, $\underline{\psi} = \bar{\underline{f}}_0^{(6)} \oplus \underline{f}_0^{(6)} = \overline{\underline{UV}}_0^{(6)} \oplus \underline{UV}_0^{(6)}$ has curvature $K = \dfrac{2}{m} = \dfrac{1}{3}$, where

$$\begin{aligned} f_0^{(6)} = [&(1+z^6, \sqrt{-1}(1-z^6), \sqrt{6}\,z(1+z^4), \sqrt{-6}\,z(1-z^4), \sqrt{15}\,z^2(-1+z^2), \\ &\sqrt{-15}\,z^2(1+z^2), \sqrt{20}\,(\frac{\sqrt{30}}{5} + \frac{\sqrt{-1}}{\sqrt{15}})z^3, \sqrt{20}\,(\frac{\sqrt{15}}{5} - \frac{\sqrt{-30}}{15})z^3)^{\mathrm{T}}] \end{aligned} \qquad .(0.30)$$

And the second fundamental form is

$$\|B\|^2 = \frac{10}{3} - \frac{72(z\bar{z})^4}{(1+z\bar{z})^8}. \qquad (0.31)$$

(2) $m = 5, K = \dfrac{2}{5}$.

Similarly, by using $V_0^{(5)} = (1, \sqrt{5}\,z, \sqrt{10}\,z^2, \sqrt{10}\,z^3, \sqrt{5}\,z^4, z^5, 0, 0)^{\mathrm{T}}$, we have the type of $W \in U(8)$ as follows:

$$W = \begin{pmatrix} 0 & 0 & 0 & 0 & 0 & -5w_{44} - 10w_{23} & 0 & 0 \\ & 0 & 0 & -\frac{1}{\sqrt{2}}w_{22} & w_{14} & 0 & w_{16} & w_{17} \\ & & w_{22} & w_{23} & -\frac{1}{\sqrt{2}}w_{33} & 0 & w_{26} & w_{27} \\ & & & w_{33} & 0 & 0 & w_{36} & w_{37} \\ & & & & 0 & 0 & w_{46} & w_{47} \\ & & * & & & 0 & 0 & 0 \\ & & & & & & w_{66} & w_{67} \\ & & & & & & & w_{77} \end{pmatrix}, \qquad (0.32)$$

and the terms in (0.17) are not all zero. We have chosen one example:

$$W = \begin{pmatrix} 0 & 0 & 0 & 0 & 0 & 1 & 0 & 0 \\ 0 & 0 & 0 & 0 & 1 & 0 & 0 & 0 \\ 0 & 0 & 0 & -\frac{3}{5} & 0 & 0 & 0 & \frac{4}{5} \\ 0 & 0 & -\frac{3}{5} & 0 & 0 & 0 & \frac{4}{5} & 0 \\ 0 & 1 & 0 & 0 & 0 & 0 & 0 & 0 \\ 1 & 0 & 0 & 0 & 0 & 0 & 0 & 0 \\ 0 & 0 & 0 & \frac{4}{5} & 0 & 0 & 0 & \frac{3}{5} \\ 0 & 0 & \frac{4}{5} & 0 & 0 & 0 & \frac{3}{5} & 0 \end{pmatrix},$$

$$U=\begin{pmatrix} \frac{1}{\sqrt{2}} & 0 & 0 & 0 & 0 & \frac{1}{\sqrt{2}} & 0 & 0 \\ \frac{\sqrt{-1}}{\sqrt{2}} & 0 & 0 & 0 & 0 & -\frac{\sqrt{-1}}{\sqrt{2}} & 0 & 0 \\ 0 & \frac{1}{\sqrt{2}} & 0 & 0 & \frac{1}{\sqrt{2}} & 0 & 0 & 0 \\ 0 & \frac{\sqrt{-1}}{\sqrt{2}} & 0 & 0 & -\frac{\sqrt{-1}}{\sqrt{2}} & 0 & 0 & 0 \\ 0 & 0 & -\frac{1}{\sqrt{2}} & -\frac{3}{5\sqrt{2}} & 0 & 0 & 0 & \frac{2\sqrt{2}}{5} \\ 0 & 0 & \frac{\sqrt{-1}}{\sqrt{2}} & \frac{3\sqrt{-1}}{5\sqrt{2}} & 0 & 0 & 0 & -\frac{2\sqrt{-2}}{5} \\ 0 & 0 & 0 & \frac{2\sqrt{2}}{5} & 0 & 0 & \frac{1}{\sqrt{2}} & \frac{3}{5\sqrt{2}} \\ 0 & 0 & 0 & \frac{2\sqrt{-2}}{5} & 0 & 0 & -\frac{\sqrt{-1}}{\sqrt{2}} & \frac{3\sqrt{-1}}{5\sqrt{2}} \end{pmatrix},$$

In this case, $\underline{\psi}=\underline{\bar{f}}_0^{(5)}\oplus\underline{f}_0^{(5)}=\underline{\overline{UV}}_0^{(5)}\oplus\underline{UV}_0^{(5)}$ has curvature $K=\frac{2}{m}=\frac{2}{5}$, where

$$f_0^{(5)}=[(1+z^5,\sqrt{-1}(1-z^5),\sqrt{5}(1+z^3),\sqrt{-5}\,z(1-z^3),\sqrt{10}\,z^2(1-\frac{3}{5}z), \\ \sqrt{-10}\,z^2(1+\frac{3}{5}z),\frac{4\sqrt{10}}{5}z^3,\frac{4\sqrt{-10}}{5}z^3)^{\mathrm{T}}]. \tag{0.33}$$

The second fundamental form $B$ of $\psi$ can be obtained:

$$\| B \|^2=\frac{16}{5}-\frac{8(z\bar{z})^3}{25(1+z\bar{z})^6}.$$

(4) $m=4, K=\frac{1}{2}$.

Similar to (1), by using $V_0^{(4)}=(1,2z,\sqrt{6}\,z^2,2z^3,z^4,0,0,0)^{\mathrm{T}}$, we have the type of $W\in U(8)$ as follows:

$$W=\begin{pmatrix} 0 & 0 & 0 & 0 & -4w_{13}-3w_{22} & 0 & 0 & 0 \\ & 0 & 0 & w_{13} & 0 & w_{15} & w_{16} & w_{17} \\ & & w_{22} & 0 & 0 & w_{25} & w_{26} & w_{27} \\ & & & 0 & 0 & w_{35} & w_{36} & w_{37} \\ & & & & 0 & 0 & 0 & 0 \\ & & * & & & w_{55} & w_{56} & w_{57} \\ & & & & & & w_{66} & w_{67} \\ & & & & & & & w_{77} \end{pmatrix}, \tag{0.34}$$

and terms in (0.18) are not all zero. An example have been chosen as follows:

$$
W=\begin{pmatrix}
0 & 0 & 0 & 0 & 1 & 0 & 0 & 0\\
0 & 0 & 0 & -\frac{1}{2} & 0 & 0 & 0 & \frac{\sqrt{3}}{2}\\
0 & 0 & \frac{1}{3} & 0 & 0 & \frac{\sqrt{8}}{3} & 0 & 0\\
0 & \frac{1}{2} & 0 & 0 & 0 & 0 & \frac{\sqrt{3}}{2} & 0\\
1 & 0 & 0 & 0 & 0 & 0 & 0 & 0\\
0 & 0 & \frac{\sqrt{8}}{3} & 0 & 0 & -\frac{1}{3} & 0 & 0\\
0 & 0 & 0 & \frac{\sqrt{3}}{2} & 0 & 0 & 0 & \frac{1}{2}\\
0 & \frac{\sqrt{3}}{2} & 0 & 0 & 0 & 0 & \frac{1}{2} & 0
\end{pmatrix},
$$

$$
U=\begin{pmatrix}
\frac{1}{\sqrt{2}} & 0 & 0 & 0 & \frac{1}{\sqrt{2}} & 0 & 0 & 0\\
-\frac{\sqrt{-1}}{\sqrt{2}} & 0 & 0 & 0 & \frac{\sqrt{-1}}{\sqrt{2}} & 0 & 0 & 0\\
0 & 0 & u_{22} & 0 & 0 & u_{25} & 0 & 0\\
0 & 0 & u_{32} & 0 & 0 & u_{35} & 0 & 0\\
0 & \frac{1}{\sqrt{2}} & 0 & -\frac{\sqrt{2}}{4} & 0 & 0 & 0 & \frac{\sqrt{6}}{4}\\
0 & \frac{\sqrt{-1}}{\sqrt{2}} & 0 & \frac{\sqrt{-2}}{4} & 0 & 0 & 0 & -\frac{\sqrt{-6}}{4}\\
0 & 0 & 0 & \frac{\sqrt{6}}{4} & 0 & 0 & \frac{1}{\sqrt{2}} & \frac{\sqrt{2}}{4}\\
0 & 0 & 0 & \frac{\sqrt{-6}}{4} & 0 & 0 & \frac{\sqrt{-1}}{\sqrt{2}} & -\frac{\sqrt{-2}}{4}
\end{pmatrix},
$$

where $\begin{cases} u_{22}=\dfrac{\sqrt{2}+\sqrt{-2}}{3} \\ u_{32}=\dfrac{2-\sqrt{-1}}{3} \end{cases}$ , and $\begin{cases} u_{25}=\dfrac{1-2\sqrt{-1}}{3} \\ u_{35}=\dfrac{\sqrt{2}+\sqrt{-2}}{3} \end{cases}$ .

In this case, $\underline{\psi}=\underline{\bar{f}}_0^{(4)}\oplus\underline{f}_0^{(4)}=\underline{\overline{UV}}_0^{(4)}\oplus\underline{UV}_0^{(4)}$ has curvature $K=\frac{1}{2}$, where

$$
\begin{aligned}
f_0^{(4)}=[&(1+z^4,\sqrt{-1}(-1+z^4),\sqrt{6}\,z^2(\frac{2}{3}+\frac{\sqrt{-2}}{3}),\sqrt{6}\,z^2(\frac{2\sqrt{2}}{3}-\frac{\sqrt{-2}}{3}),\\
&2z-z^3,\sqrt{-1}(2z+z^3),\sqrt{3}\,z^3,\sqrt{-3}\,z^3)^{\mathrm{T}}].
\end{aligned}
\tag{0.35}
$$

And the second fundamental form is

$$
\|B\|^2=3-\frac{2(z\bar{z})^2}{(1+z\bar{z})^4}.
$$

(3) $m=3, K=\dfrac{2}{3}$.

Similar to (1) and $V_0^{(3)}=(1,\sqrt{3}\,z,\sqrt{3}\,z^2,z^3,0,0,0,0)^{\mathrm{T}}$, we have

$$W=\begin{pmatrix} 0 & 0 & 0 & -3w_{12} & 0 & 0 & 0 & 0 \\ & 0 & w_{12} & 0 & w_{14} & w_{15} & w_{16} & w_{17} \\ & & 0 & 0 & w_{24} & w_{25} & w_{26} & w_{27} \\ & & & 0 & 0 & 0 & 0 & 0 \\ & & & & w_{44} & w_{45} & w_{46} & w_{47} \\ & & * & & & w_{55} & w_{56} & w_{57} \\ & & & & & & w_{66} & w_{67} \\ & & & & & & & w_{77} \end{pmatrix}, \tag{0.36}$$

and the terms in (0.19) are not all zero. An example have been chosen as follows:

$$W=\begin{pmatrix} 0 & 0 & 0 & -\frac{3}{5} & 0 & 0 & 0 & \frac{4}{5} \\ 0 & 0 & \frac{1}{5} & 0 & 0 & 0 & -\frac{2\sqrt{6}}{5} & 0 \\ 0 & \frac{1}{5} & 0 & 0 & 0 & \frac{2\sqrt{6}}{5} & 0 & 0 \\ -\frac{3}{5} & 0 & 0 & 0 & \frac{4}{5} & 0 & 0 & 0 \\ 0 & 0 & 0 & \frac{4}{5} & 0 & 0 & 0 & \frac{3}{5} \\ 0 & 0 & \frac{2\sqrt{6}}{5} & 0 & 0 & 0 & \frac{1}{5} & 0 \\ 0 & -\frac{2\sqrt{6}}{5} & 0 & 0 & 0 & \frac{1}{5} & 0 & 0 \\ \frac{4}{5} & 0 & 0 & 0 & \frac{3}{5} & 0 & 0 & 0 \end{pmatrix},$$

$$U=\begin{pmatrix}
\frac{1}{\sqrt{2}} & 0 & 0 & -\frac{3\sqrt{2}}{10} & 0 & 0 & 0 & \frac{2\sqrt{2}}{5} \\
-\frac{\sqrt{-1}}{\sqrt{2}} & 0 & 0 & \frac{3\sqrt{-2}}{10} & 0 & 0 & 0 & -\frac{2\sqrt{-2}}{5} \\
0 & 0 & 0 & \frac{2\sqrt{2}}{5} & \frac{1}{\sqrt{2}} & 0 & 0 & \frac{3\sqrt{2}}{10} \\
0 & 0 & 0 & \frac{2\sqrt{-2}}{5} & \frac{\sqrt{-1}}{\sqrt{2}} & 0 & 0 & \frac{3\sqrt{-2}}{10} \\
0 & \frac{1}{\sqrt{2}} & \frac{\sqrt{2}}{10} & 0 & 0 & 0 & -\frac{2\sqrt{3}}{5} & 0 \\
0 & \frac{\sqrt{-1}}{\sqrt{2}} & -\frac{\sqrt{-2}}{10}0 & 0 & 0 & 0 & \frac{2\sqrt{-3}}{5} & 0 \\
0 & 0 & -\frac{2\sqrt{-3}}{5} & 0 & 0 & \frac{\sqrt{-1}}{\sqrt{2}} & -\frac{2\sqrt{-2}}{10} & 0 \\
0 & 0 & -\frac{2\sqrt{3}}{5}0 & 0 & 0 & \frac{\sqrt{-1}}{\sqrt{2}} & -\frac{\sqrt{2}}{10} & 0
\end{pmatrix}.$$

In this case, $\underline{\psi}=\underline{\bar{f}}_0^{(3)}\oplus\underline{f}_0^{(3)}=\underline{\overline{UV}}_0^{(3)}\oplus\underline{UV}_0^{(3)}$ has curvature $K=\frac{2}{3}$, where

$$\begin{aligned}f_0^{(3)}=[(1-\frac{3}{5}z^3,\sqrt{-1}(1+\frac{3}{5}z^3),\frac{4}{5}Z^3,\frac{4\sqrt{-1}}{5}z^3,\sqrt{3}z(1+\frac{z}{5}),\\ \sqrt{-3}z(1-\frac{z}{5}),-\frac{2\sqrt{-6}}{5}z^2,-\frac{2\sqrt{6}}{5}z^2)^{\mathrm{T}}].\end{aligned} \tag{0.37}$$

Direct computations give

$$\|B\|^2=\frac{8}{3}-\frac{16z\bar{z}}{25(1+z\bar{z})^2}.$$

(5) $m=2,K=1$.

Since $V_0^{(2)}=(1,\sqrt{2}z,z^2,0,0,0,0,0)^{\mathrm{T}}$, we obtain that $f_0^{(2)}=UV_0^{(2)}$ is not linearly full in $Q_6$ for all $U\in U(8)$, and therefore, this case does not exist.

In summary, we have arrived at the conclusion.

## 3.2 finite isotropy order $r=3$

**Lemma 3.3** *Let* $\psi$ *is a linearly full reducible harmonic map with constant curvature* $K$ *of isotropy order* $r=3$. *In fact,* $\underline{\psi}=\underline{\overline{f_0^{(m)}}}\oplus\underline{f_0^{(m)}}(4\le m\le 7)$, *where* $f_0^{(m)}:S^2\to\mathbb{C}P^{(m)}$ *is holomorphic. Then* $m\neq 7$ *and* $m\neq 6$.

Proof: Through the previous argument, we can know such $\psi$ satisfying (3.2). it is easy to verify the fact that (0.12) is equivalent to

$$\begin{cases} \langle \partial^{(i)} f_0^{(m)}, \bar{f}_0^{(m)} \rangle = 0, \quad (0 \le i \le 3) \\ \langle \partial^{(4)} f_0^{(m)}, \bar{f}_0^{(m)} \rangle \neq 0. \end{cases}$$

Then with some simple computations, we have $\langle \partial^{(3)} f_0, \bar{f}_0^{(m)} \rangle = 0$ if the formulas $\langle f_0^{(m)}, \bar{f}_0^{(m)} \rangle = 0$ and $\langle \partial^{(2)} f_0^{(m)}, \bar{f}_0^{(m)} \rangle = 0$ hold true. From Lemma 3.2, we have already obtained that $\langle f_0^{(m)}, \bar{f}_0^{(m)} \rangle = 0$ and $\langle f_1^{(m)}, \bar{f}_0^{(m)} \rangle = 0$ are equivalent. So the relations above are converted into the following formulas:

$$\begin{cases} \langle f_0^{(m)}, \overline{f_i^{(m)}} \rangle = 0, \quad (i = 0, 2) \\ \langle f_0^{(m)}, \overline{f_4^{(m)}} \rangle \neq 0. \end{cases} \tag{0.38}$$

When $m = 7$ and $K = 2/7$, similarly to Lemma 3.2, by the expression of $V_0^{(7)}$, $\partial^{(2)} V_0^{(7)}$ and $\partial^{(4)} V_0^{(7)}$, using the method of indeterminate coefficients, we can have

$$\begin{gathered} w_{04} = \frac{\sqrt{35}}{10} w_{22}, \quad w_{05} = \frac{\sqrt{35}}{2} w_{23}, \quad w_{13} = -\frac{2}{\sqrt{5}} w_{22}, \quad w_{14} = -\frac{3\sqrt{3}}{2} w_{23}, \\ w_{06} = \frac{3\sqrt{105}}{5} w_{24} + 2\sqrt{7}\, w_{33}, \quad w_{07} = 14 w_{25} + 35 w_{34}, \\ w_{15} = -\frac{8}{\sqrt{5}} w_{24} - \frac{3\sqrt{3}}{2} w_{33}, \quad w_{16} = -5 w_{25} - 10 w_{34}, \\ w_{17} = \frac{3\sqrt{105}}{5} w_{35} + 2\sqrt{7}\, w_{44}, \quad w_{26} = -\frac{8}{\sqrt{5}} w_{35} - \frac{3\sqrt{3}}{2} w_{44}, \\ w_{37} = \frac{\sqrt{35}}{10} w_{55}, \quad w_{27} = \frac{\sqrt{35}}{2} w_{45}, \quad w_{46} = -\frac{2}{\sqrt{5}} w_{55}, \quad w_{36} = -\frac{3\sqrt{3}}{2} w_{45}. \end{gathered} \tag{0.39}$$

Using (0.39) and the property of unitarity, it is straightforward computation to derive

$$W = \begin{pmatrix} 0 & 0 & 0 & 0 & 0 & 0 & 0 & 14w_{25} + 35w_{34} \\ & 0 & 0 & 0 & 0 & 0 & -5w_{25} - 10w_{34} & 0 \\ & & 0 & 0 & 0 & w_{25} & 0 & 0 \\ & & & 0 & w_{34} & 0 & 0 & 0 \\ & & & & 0 & 0 & 0 & 0 \\ & & * & & & 0 & 0 & 0 \\ & & & & & & 0 & 0 \\ & & & & & & & 0 \end{pmatrix},$$

with $|w_{25}| = |w_{34}| = 1$, which is contradicts $|14 w_{25} + 35 w_{34}| = 1$. So $m \neq 7$ is proved.

Analogously, $m \neq 6$ can be proved by using the same method above.

Let

$$W = \begin{pmatrix} & & & & w_{04} & w_{05} & w_{06} & w_{07} \\ & 0 & & & w_{14} & w_{15} & w_{16} & w_{17} \\ & & & & w_{24} & w_{25} & w_{26} & w_{27} \\ & & & & w_{34} & w_{35} & w_{36} & w_{37} \\ w_{40} & w_{41} & w_{42} & w_{43} & w_{44} & w_{45} & w_{46} & w_{47} \\ w_{50} & w_{51} & w_{52} & w_{53} & w_{54} & w_{55} & w_{56} & w_{57} \\ w_{60} & w_{61} & w_{62} & w_{63} & w_{64} & w_{65} & w_{66} & w_{67} \\ w_{70} & w_{71} & w_{72} & w_{73} & w_{74} & w_{75} & w_{76} & w_{77} \end{pmatrix} \tag{0.40}$$

By Lemma 3.2 , Lemma 3.3 and [10, Proposition 3.5], we conclude a rough classification theorem of linearly full conformal minimal immersion with constant curvature from $S^2$ to $G(2,8;\mathbb{R})$ as follows.

**Theorem 1** *Let $\psi: S^2 \to G(2,8;\mathbb{R})$ be a linearly full reducible conformal minimal immersion with constant curvature $K$ of isotropy order $r$. Then, $r=1$ or $r=3$ or $r=\infty$, up to an isometry of $G(2,8;\mathbb{R})$,*

*(1) when $r=1$, see Lemma 3.2;*

*(2) when $r=3$, see Lemma 3.3;*

*(3) when $r=\infty$:*

*a) $\underline{\psi} = \underline{UV_3^{(6)}} \oplus \underline{c_0}$ with $K=\frac{1}{6}$ and $\|B\|^2=\frac{5}{3}$ for some $U \in U(7)$ and $c_0=(0,0,0,0,0,0,1)^{\mathrm{T}}$.*

*b) $\underline{\psi} = \overline{\underline{UV_0^{(3)}}} \oplus \underline{UV_0^{(3)}}$ with $K=\frac{2}{3}$ and $\|B\|^2=\frac{8}{3}$ for some $U \in G_W$, where $W$ has the form (0.40).*

*In cases of (1), b), there are many different types of $W$, and thus there exist different $U \in U(8)$ such that $UV_0^{(m)}$ are linearly full in $G(2,8;\mathbb{R})$, and they are not $SO(8)$-equivalent.*

Here, the type of $U$ in a) and b) of Theorem 1 can be chosen to be the forms, respectively.

$$U_1 = \begin{pmatrix} \frac{1}{\sqrt{2}} & 0 & 0 & 0 & 0 & 0 & \frac{1}{\sqrt{2}} \\ \frac{\sqrt{-1}}{\sqrt{2}} & 0 & 0 & 0 & 0 & 0 & -\frac{\sqrt{-1}}{\sqrt{2}} \\ 0 & \frac{1}{\sqrt{2}} & 0 & 0 & 0 & -\frac{1}{\sqrt{2}} & 0 \\ 0 & \frac{\sqrt{-1}}{\sqrt{2}} & 0 & 0 & 0 & \frac{\sqrt{-1}}{\sqrt{2}} & 0 \\ 0 & 0 & \frac{1}{\sqrt{2}} & 0 & \frac{1}{\sqrt{2}} & 0 & 0 \\ 0 & 0 & \frac{\sqrt{-1}}{\sqrt{2}} & 0 & -\frac{\sqrt{-1}}{\sqrt{2}} & 0 & 0 \\ 0 & 0 & 0 & \sqrt{-1} & 0 & 0 & 0 \end{pmatrix}, \tag{0.41}$$

$$U_2 = \begin{pmatrix} \frac{1}{\sqrt{2}} & 0 & 0 & 0 & \\ \frac{\sqrt{-1}}{\sqrt{2}} & 0 & 0 & & \\ 0 & \frac{1}{\sqrt{2}} & 0 & 0 & \\ 0 & \frac{\sqrt{-1}}{\sqrt{2}} & 0 & 0 & \\ 0 & 0 & \frac{1}{\sqrt{2}} & 0 & * \\ 0 & 0 & \frac{\sqrt{-1}}{\sqrt{2}} & 0 & \\ 0 & 0 & 0 & \frac{1}{\sqrt{2}} & \\ 0 & 0 & 0 & \frac{\sqrt{-1}}{\sqrt{2}} & \end{pmatrix}. \tag{0.42}$$

Then up to an isometry of $G(2,8;\mathbb{R})$, either

$$\underline{\psi} = \underline{UV_3^{(6)}} \oplus \underline{c_0} = \underline{f_3^{(6)}} \oplus \underline{c_0}$$

with

$$\begin{aligned} f_3^{(6)} = [(&\sqrt{10}(z^3 - \bar{z}^3), \sqrt{-10}(-z^3 - \bar{z}^3), \sqrt{15}(\bar{z}^2 - z^2)(1 - |z|^2), \\ &\sqrt{-15}(\bar{z}^2 + z^2)(1 - |z|^2), \sqrt{6}(z - \bar{z})(1 - 3|z|^2 + |z|^4), \\ &\sqrt{-6}(z + \bar{z})(-1 + 3|z|^2 - |z|^4), \sqrt{-1}(1 - 9|z|^2 + 9|z|^4 - |z|^6)^{\mathrm{T}}] \end{aligned} \tag{0.43}$$

or

$$\underline{\psi} = \overline{\underline{UV_0^{(3)}}} \oplus \underline{UV_0^{(3)}} = \overline{\underline{f_0^{(3)}}} \oplus \underline{f_0^{(3)}}$$

with

$$f_0^{(3)} = [1, \sqrt{-1}, \sqrt{3}z, \sqrt{-3}z.\sqrt{3}z^2, \sqrt{-3}z^2, z^3, \sqrt{-1}z^3)^{\mathrm{T}}]. \tag{0.44}$$